\documentclass[12pt, a4paper]{amsart}
\usepackage{amscd,amsthm,amsfonts,amssymb,amsmath,euscript}
\usepackage[dvips]{graphicx}
\usepackage[dvips]{graphics}
\usepackage[matrix,arrow]{xy}
\usepackage{longtable}
\usepackage{color}
\usepackage{url}

\newcounter{statements}

\newtheorem{theorem}[statements]{Theorem}
\newtheorem{conjecture}[statements]{Conjecture}
\newtheorem{problem}[statements]{Problem}

\theoremstyle{definition}

\newtheorem{definition}[statements]{Definition}

\theoremstyle{remark}

\newtheorem{remark}[statements]{Remark}

\newcommand{\ZZ}{{\mathbb Z}}

\newcommand{\PP}{{\mathbb P}}
\newcommand{\CC}{{\mathbb C}}
\newcommand{\RR}{{\mathbb R}}

\newcommand{\TT}{{\mathbb T}}

\newcommand{\Aff}{{\mathbb A}}

\newfont{\smallskob}{cmbx7 scaled\magstep4}
\newfont{\bigskob}{cmbx12 scaled\magstep4}

\newcommand{\Pic}{\mathrm{Pic}\,}

\newcommand{\iX}{i_X}

\begin{document}


\title{On Calabi--Yau compactifications of toric Landau--Ginzburg models for Fano complete intersections}

\author{Victor Przyjalkowski}

\thanks{This work is supported by the Russian Science Foundation under grant 14-50-00005.}

\address{Steklov Mathematical Institute of Russian Academy of Sciences, 8 Gubkina street, Moscow 119991, Russia} %

\email{victorprz@mi.ras.ru, victorprz@gmail.com}


\maketitle


\begin{abstract}
Toric Landau--Ginzburg models of Givental's type for Fano complete intersections are known to have Calabi--Yau compactifications.
We give an alternative proof of this fact. As an output of our proof
we get a description of fibers over infinity for compactified
toric Landau--Ginzburg models.
\end{abstract}




\bigskip

\section{Introduction}

In~\cite{Gi97} (see also~\cite{HV00}) Givental presented Landau--Ginzburg models for complete intersections with non-positive canonical classes in smooth toric varieties.
They are certain linear sections in complex tori equipped by complex-valued functions called \emph{superpotentials}.
He proved that the regularized $I$-series, that are generating series for genus 0 one-pointed Gromov--Witten invariants with descendants, are solutions
of Picard--Fuchs equations for families of fibers of superpotentials (see more details below). In some cases, such as for
Fano complete intersections in projective spaces, Givental's models can be birationally presented as complex tori with
complex valued functions on them. These functions are (in some bases) given by Laurent polynomials. Such Laurent polynomials are
\emph{toric Landau--Ginzburg models}, which, besides the $I$-series-periods correspondence (see~\cite{Gi97} and~\cite{Prz13}),
means that there is a correspondence between Laurent polynomials and toric degenerations of the complete intersections
(see~\cite{ILP13}), and that they admit \emph{Calabi--Yau compactifications} (see~\cite{Prz13}), see details in Definition~\ref{definition: toric LG}.
(These Calabi--Yau compactifications are described in details in~\cite{PSh15a}.)

More precise, a Landau--Ginzburg model from Homological Mirror Symmetry point of view is a smooth quasiprojective variety with
a complex-valued function (superpotential). This variety and the superpotential satisfy
a number of properties, algebraic and symplectic. In particular, fibers of the superpotential should be
Calabi--Yau varieties. Compactification Principle (see~\cite[Principle 32]{Prz13}) states that
``correct'' toric Landau--Ginzburg model, considered as a family of fibers of the superpotential, admits a fiberwise compactification to
a Landau--Ginzburg model satisfying Homological Mirror Symmetry conjecture. In particular this compactification is
a family of compact Calabi--Yau varieties with smooth total space.

A Landau--Ginzburg model by definition is a family of varieties over $\Aff^1$.
However sometimes a compactification to a family over $\PP^1$ is important,
see, for instance,~\cite{AKO06},~\cite{KKP17},~\cite{LP16}.
In~\cite{Prz16} such compactifications of toric Landau--Ginzburg models
for smooth Fano threefolds are studied in a systematic way. That is,
for any Fano threefold one can choose a ``good'' toric Landau--Ginzburg model $f$ such that the following holds.
Let $\Delta$ be a Newton polytope for $f$ and let $\nabla$ be a polytope dual to $\Delta$.
The polytope $\Delta$ is reflexive, that is $\nabla$ is integral, and thus it defines a toric Fano variety $T$. Fibers of the map determined by $f$ are given
by Laurent polynomials supported in $\Delta$. This means that they correspond to elements of an anticanonical
linear system of $T$. In particular they have trivial canonical classes. One can check that
$T$ has a crepant (toric) resolution $\widetilde T$, and a base locus of the family
of (compactifications in $\widetilde T$ of) fibers of $f$ is a union of smooth (rational) curves, possibly with multiplicities.
(This is provided by the fact that $f$ is \emph{Minkowski Laurent polynomial}, see~\cite{CCGK16}.)
 This means that a resolution of the base locus is the required compactification.
In addition to a Calabi--Yau compactification over $\Aff^1$, this procedure gives a description of the fiber over infinity as a boundary divisor
of $\widetilde T$.
It is reduced, which means that the mo\-no\-dro\-my of the compactified Landau--Ginzburg model at infinity is maximally unipotent by Griffiths--Landman--Grothendieck Theorem. In other words, the compactifications are \emph{tame compactified Landau--Ginzburg models}
(see~\cite{KKP17}).

In the paper we follow the method of~\cite{Prz16} to reproduce Calabi--Yau compactifications for Givental's
Landau--Ginzburg models for complete intersections constructed in~\cite{Prz13} and~\cite{PSh15a}. Newton polytopes of the toric Landau--Ginzburg models are reflexive in our case
(see~\cite{ILP13} and below), and the base loci of (compactifications in dual toric Fano varieties of) the families under consideration
are unions of (rational) varieties. In the opposite to the threefold case components of the base loci can be singular.
However they become smooth after crepant resolutions of the toric Fano varieties.
This is provided by symmetric ``binomial'' coefficients of the toric Landau--Ginzburg models.
As a result of this construction one gets log Calabi--Yau compactifications (over $\PP^1$), which give the required Calabi--Yau compactifications over $\Aff^1$.

A side effect of these compactifications is a precise description of fibers over infinity.
That is, one can define a certain $n$-dimensional polytope, whose vertices have coordinates given by degrees $d_1,\ldots,d_k$ of hypersurfaces that cut
a complete intersection of dimension $n$. A fiber over infinity of a compactified toric Landau--Ginzburg model for the complete intersection can be described by a maximal integral triangulation of the polytope.
Vertices of simplices of the triangulation are integral points on the boundary of the polytope. These vertices correspond to components of the fiber over infinity,
edges of the triangulation correspond to intersections of the components, etc. In particular, fibers over infinity for
different compactifications (obtained in a way described above) differ by flops of a total space.
The compactifications are tame compactified Landau--Ginzburg models, so one can apply results from~\cite{KKP17} for them.

\begin{theorem}
\label{theorem:main}
Let $X\subset \PP^N$ be a Fano complete intersection of hypersurfaces of degrees $d_1,\ldots,d_k$.
Let $\iX=N+1-\sum d_i$. Let $f_X$ be a toric Landau--Ginzburg model of Givental's type for $X$.
Then $f_X$ admits a log Calabi--Yau compactification $f_X\colon Z\to \PP^1$ such that $f_X^{-1}(\infty)$
is a reduced divisor, which is a union of smooth rational varieties.
It consists of $k_{d_1,\ldots,d_k;\iX}$ components (see Definition~\ref{definition: number of components})
and combinatorially it is given by a triangulation of a sphere.
\end{theorem}

\section{Preliminaries}
First let us give a definition of a toric Landau--Ginzburg model. More details see, say, in~\cite{Prz13}.
Let $X$ be a smooth Fano variety of dimension $n\geq 3$ with $\Pic (X)=\ZZ$.
Let $\mathbf 1$ be the fundamental class of $X$.
The series
$$
I^{X}_{0}(t)
=1+\sum_{d>0,\ a\in \ZZ_{\geq 0}} d! \langle\tau_{a} \mathbf 1\rangle_{d}
\cdot t^{d},
$$
where $\langle\tau_{a} \mathbf 1\rangle_{d}$ is a \emph{one-pointed genus $0$ Gromov--Witten invariant with descendants} for anticanonical
degree $d$ curves on $X$, see~\cite[VI-2.1]{Ma99},
is called \emph{a constant term of regularized $I$-series} (or \emph{a constant term of regularized Givental's $J$-series}) for $X$.

\begin{definition}
Let $f$ be a Laurent polynomial in $n$ variables $x_1,\ldots,x_n$.
The integral
$$
I_f(t)=
\frac{1}{(2\pi i)^n}\int\limits_{|x_i|=\varepsilon_i}\frac{dx_1}{x_1}\wedge\ldots\wedge \frac{dx_n}{x_n}\frac{1
}{1-tf} = \frac{1}{(2\pi i)^n}\sum_{j=0}^\infty  t^j \cdot \int\limits_{|x_i|=\varepsilon_i}
f^j\frac{dx_1}{x_1}\wedge\ldots\wedge \frac{dx_n}{x_n}
 \in\CC[[t]],
$$
where $\varepsilon_i$ are arbitrary positive real numbers, is called \emph{the main period} for $f$.
\end{definition}

\begin{remark}
\label{remark: Picard--Fuchs}
Let~$\phi[f]$ be a constant term of a Laurent polynomial~$f$.
Then $I_f(t)=\sum \phi[f^j] t^j$.
\end{remark}

The following theorem (which is a mathematical folklore, see~\cite[Proposition 2.3]{Prz08} 
for the proof)
justifies this definition.

\begin{theorem}
\label{theorem: Picard--Fuchs}
Let $f$ be a Laurent polynomial in $n$ variables.
Let $P$ be a Picard--Fuchs differential operator for a pencil of hypersurfaces in a torus
provided by $f$.
Then
one has~\mbox{$P[I_f(t)]=0$}.
\end{theorem}

For a Laurent polynomial $f\in \CC[\ZZ^n]$ we denote its Newton polytope, i.\,e. a
convex hull in $\ZZ^n\otimes \RR$ of non-zero monomials of $f$, by $N(f)$.
Given a toric variety $T$ we define
\emph{a fan} (or \emph{spanning}) polytope $F(T)$ as a convex hull of integral generators
of fan's rays for $T$.

\begin{definition}[{see~\cite[\S6]{Prz13}}]
\label{definition: toric LG}
\emph{A toric Landau--Ginzburg model} for a smooth Fano variety $X$ of dimension $n$ (corresponding to the anticanonical class) is a Laurent polynomial $\mbox{$f\in \TT[x_1, \ldots, x_n]$}$ which satisfies the following.
\begin{description}
  \item[Period condition] One has $I_f(t)=\widetilde{I}_0^{X}(t)$.
  \item[Calabi--Yau condition] There exists a relative compactification of a family
$$f\colon (\CC^*)^n\to \CC,$$
whose total space is a (non-compact) smooth Calabi--Yau
variety $Y$. Such compactification is called \emph{a Calabi--Yau compactification}.
  \item[Toric condition] There is a degeneration
   $X\rightsquigarrow T$ to a toric variety~$T$ such that $F(T)=N(f)$.
\end{description}
\end{definition}

\begin{definition}
A compactification of the family $f\colon (\CC^*)^n\to \CC$ to a family $f\colon Z\to \PP^1$, where $Z$ is smooth and compact and
$-K_Z=f^{-1}(\infty)$, is called a \emph{log Calabi--Yau compactification}.
\end{definition}

\begin{conjecture}[{see~\cite[Conjecture 38]{Prz13}}]
\label{conjecture:MS}
Any smooth Fano variety has a toric Landau--Ginzburg model. 
\end{conjecture}

This conjecture holds for smooth Fano threefolds (see~\cite{Prz13} and~\cite{ILP13} for the Picard rank one case
and~\cite{fanosearch},~\cite{IKKPS},~\cite{DHKLP}, and~\cite{Prz16} for the general case) and complete
intersections (see~\cite{Prz13},~\cite{ILP13}, and also~\cite{PSh15a}).
The existence of Laurent polynomials satisfying
the period condition for smooth toric varieties is shown in~\cite{Gi97},
for complete intersections in Grassmannians it is shown in~\cite{PSh14},~\cite{PSh14b}~\cite{PSh15b},
for some complete intersections in some toric varieties it is shown in~\cite{CKP14} and~\cite{DH15}.
In~\cite{DH15} the toric condition for some complete intersections in toric varieties and partial
flag manifolds is also checked.

\begin{definition}
Let us fix natural numbers $d_1,\ldots,d_k$, and $N$ such that the number $I=N+1-\sum d_j$ is positive.
The Laurent polynomial
$$
f_{d_1,\ldots,d_k;I}=\frac{\prod_{i=1}^k(x_{i,1}+\ldots+x_{i,d_i-1}+1)^{d_i}}{\prod_{i=1}^k \prod_{j=1}^{d_i-1} x_{i,j}\prod_{j=1}^{I-1} y_j}+y_1+\ldots+y_{I-1}
$$
is called \emph{standard of type $(d_1,\ldots,d_k;I)$}.
\end{definition}

Let $X\subset\PP^N$ be a Fano complete intersection of hypersurfaces of degrees $d_1,\ldots,d_k$ and let $\iX=N+1-\sum d_j$ be its index.
Its toric Landau--Ginzburg model is $f_X=f_{d_1,\ldots,d_k;\iX}$, see~\cite[\S 3.2]{Prz13}.
A toric Fano variety $T_X$ whose fan polytope is $\Delta=N(f_X)$ is called \emph{associated with $f$};
its precise description one can see in~\cite{ILP13}.
Let
$$
\nabla=\{x\ |\ \langle x,y\rangle \geq -1 \mbox{ for all } y\in \Delta \}\subset M_\RR=N^\vee\otimes \RR
$$
be the dual to $\Delta$ polytope. One gets that $\nabla$ is integral from the precise description of $\Delta$; 
coordinates of its vertices are given by rows of the matrix $M_{d_1,\ldots,d_k;\iX}$ below.
In other words,
$T_X$ and $\Delta$ are \emph{reflexive}.

Another important property of the polynomial $f_{d_1,\ldots,d_k;I}$ is related to its coefficients.
Let $Q$ be a face of $\Delta$ and let $f_{{d_1,\ldots,d_k;I},Q}$ be a restriction of $f_{d_1,\ldots,d_k;I}$ to $Q$,
which means that $f_{{d_1,\ldots,d_k;I},Q}$ is a sum of those coefficients of $f_{d_1,\ldots,d_k;I}$ that lie in $\Delta$.
Then $f_{{d_1,\ldots,d_k;I},Q}$ is standard again.

\section{Calabi--Yau compactifications}
In this section we prove Theorem~\ref{theorem:main}. We follow the compactification procedure described in~\cite{Prz16};
for details, definitions, and background see therein.
Let as above $X\subset \PP^N$ be a Fano complete intersection of hypersurfaces of degrees $d_1,\ldots,d_k$, let $\iX=N+1-\sum d_j$,
let $f_X=f_{d_1,\ldots,d_k;\iX}$, let $T_X$ be a toric Fano variety, associated with $f_X$, 
and let $\Delta=F(T_X)$.
Let $\nabla$ be a polytope dual to $\Delta$. As we mentioned above, $\nabla$ is integral (see below for the precise description of $\nabla$).
Let $T_X^\vee$ ba a toric Fano variety such that $F(T_X^\vee)=\nabla$. Anticanonical linear system on $T_X^\vee$ can
be described as (a projectivisation of) a linear system of Laurent polynomials whose Newton polytopes lie in $\Delta$.
Thus, a family $\{f_X=\lambda|\ \lambda\in \CC\}$ can be compactified to a family of (singular) Calabi--Yau anticanonical sections of $T_X^\vee$.
(We use for simplicity the same notation $f_X$ for the function on $(\CC^*)^{N-k}$ and on all its compactifications.)
It is generated by a general fiber (that is $f_X^{-1}(\lambda)$ for general $\lambda\in \CC$) and by a fiber over infinity $f_X^{-1}(\infty)$ which is nothing but a boundary
divisor $D$ of $T_X^\vee$.

To get a log Calabi--Yau compactification one needs to resolve singularities of $T_X^\vee$ and a base locus $B=f_X^{-1}(\lambda)\cap D$
of the pencil and to check that all resolutions are crepant.

Let $D=D_1\cup \ldots \cup D_r$, where $D_i$ are irreducible components of $D$. Then $B=B_1\cup\ldots \cup B_r$,
where $B_i=f_X^{-1}(\lambda) \cap D_i$. A boundary divisor $D_i$ corresponds to a vertex of $\nabla$ and, thus,
to a face $Q_i$ of $\Delta$. As a variety it can be given as follows. Let $z_1,\ldots,z_s$ be formal variables
corresponding to integral points of $Q_i$, and let $R_1,\ldots,R_p$ be homogenous binomials in $z_j$ corresponding
to relations on integral points of $Q_i$. Then $D_i=\{R_j=0|\, j=1,\ldots, p\}\subset \PP(z_1,\ldots, z_s)$.
In particular this means that $D_i$ is a cone over Segre embedded product of Veronese embedded projective spaces
or (in cases when $\iX=1$ or $\iX=2$) a Segre embedded product of Veronese embedded projective spaces itself.
We call these boundary divisors of type I and II correspondingly. The restriction of $f_X$ to $Q_i$ is a sum of
monomials (with some coefficients) that correspond to integral points of $Q_i$ and thus to variables $z_j$; the sum of $z_j$ with the coefficients of $f_X$
is a linear section of $D_i\subset \PP(z_1,\ldots, z_s)$ that gives $B_i$.

If $D_i$ is of type I, then one of variables $z_j$, say $z_s$,
corresponds to an element of the basis of the lattice $N$ that is associated to a variable $y_q$ for some $q=1,\ldots, \iX-1$.
This means that $B_i$ projects isomorphically to its image under projection $\PP(z_1,\ldots, z_s)\to \PP(z_1,\ldots, z_{s-1})$,
and the image is 
(a cone over) Segre embedded product of the same Veronese embedded projective spaces as ones for $D_i$.
This in particular means that ``singularities of $B_i$ come from $D_i$'', so after a resolution
of $T_X^\vee$ the component $B_i$ becomes smooth. Now let $D_i$ be of type II. Then $B_i$ is a union of cones of linear sections (with multiplicities)
of the projective spaces, so it is a union of smooth projective spaces (with multiplicities), cf.~\cite[Lemma~27]{Prz16}.

Consider a toric variety $\widetilde T_X^\vee$ given by a fan obtained by a triangulation of $\nabla$, vertices of whose simplices are
all integral points on the boundary of $\nabla$. Cones of this fan are generated by elements of the basis of $N$ (that ones that correspond to variables $y_i$)
and elements of (products of) standard triangulations of equilateral simplices. This means that they form a part of a basis of $N$,
so $\widetilde T_X^\vee$ is smooth and it is a crepant resolution of $T_X^\vee$. Moreover, after this resolution
$B$ becomes a union of smooth rational varieties (possibly with multiplicities). After crepant resolution $Z\to \widetilde T_X^\vee$ of this base locus (crepancy is provided by smoothness of $\widetilde T_X^\vee$ and components of $B$ and the fact that $B\in \widetilde T_X^\vee$ is of pure codimension two, see~\cite[Proposition 28]{Prz16}) one gets the required log Calabi--Yau compactification $Z$.

\begin{remark}
One can define weak Landau--Ginzburg models (that is Laurent polynomials satisfying the period condition) for smooth well formed Fano
complete intersections in weighted projective spaces provided that they admits so called \emph{nice nef partitions}.
By~\cite{ILP13} toric condition holds for them.
Such nef partitions exist for complete intersections of Cartier divisors (see~\cite{Prz11}) and complete intersections
of codimension not greater then $2$ (see~\cite{PSh17}). The proof of Theorem~\ref{theorem:main} can be repeated for
these weak Landau--Ginzburg models provided that their Newton polytopes are reflexive (which is in fact not a common case).
\end{remark}

Now let us describe polytopes $\Delta$ and $\nabla$ precisely.
Vertices of $\Delta\subset N$ are all possible vectors $u$ and $v$ such that
$$
u=\left(u_{1,1},\ldots, u_{1,d_1-1},\ldots,u_{k,1},\ldots, u_{k,d_k-1},-1,\ldots,-1\right),
$$
(so that the last $\iX-1$ coordinates of $u$ equal $-1$),
where
$$
\left(u_{i,1},\ldots, u_{i,d_i-1}\right)=
\left(-1,\ldots,d_i-1,\ldots,-1\right)
$$
with $d_i-1$ standing on all $d_i-1$ possible places, or
$$
\left(u_{i,1},\ldots, u_{i,d_i-1}\right)=
\left(-1,\ldots,-1\right),
$$
and
$v=(0,\ldots,0,1,0,\ldots,0)$,
where $1$ stands on one of the last $\iX-1$ places.

Vertices of $\nabla$ are rows of the matrix
$$
M_{d_1,\ldots,d_k;\iX}=\left(%
\begin{array}{rrrr|r|rrrr|rrr}
  \iX & 0 & \ldots & 0 & \ldots & 0 & 0 & \ldots & 0 & -1 & \ldots & -1 \\
  0 & \iX & \ldots & 0 & \ldots & 0 & 0 & \ldots & 0 & -1 & \ldots & -1 \\
  \ldots & \ldots & \ldots & \ldots & \ldots & \ldots & \ldots & \ldots & \ldots & \ldots & \ldots & \ldots\\
  0 & 0 & \ldots & \iX & \ldots & 0 & 0 & \ldots & 0 & -1 & \ldots & -1 \\
  -\iX & -\iX & \ldots & -\iX & \ldots & 0 & 0 & \ldots & 0 & -1 & \ldots & -1 \\
  \hline
  \ldots & \ldots & \ldots & \ldots & \ldots & \ldots & \ldots & \ldots & \ldots & \ldots & \ldots & \ldots\\
  \hline
  0 & 0 & \ldots & 0 & \ldots & \iX & 0 & \ldots & 0 & -1 & \ldots & -1 \\
  0 & 0 & \ldots & 0 & \ldots & 0 & \iX & \ldots & 0 & -1 & \ldots & -1 \\
  \ldots & \ldots & \ldots & \ldots & \ldots & \ldots & \ldots & \ldots & \ldots & \ldots & \ldots & \ldots\\
  0 & 0 & \ldots & 0 & \ldots & 0 & 0 & \ldots & \iX & -1 & \ldots & -1 \\
  0 & 0 & \ldots & 0 & \ldots & -\iX & -\iX & \ldots & -\iX& -1 & \ldots & -1 \\
\hline
  0 & 0 & \ldots & 0 & \ldots & 0 & 0 & \ldots & 0 & \iX-1 & \ldots & -1 \\
  \ldots & \ldots & \ldots & \ldots & \ldots & \ldots & \ldots & \ldots & \ldots & \ldots & \ldots & \ldots \\
  0 & 0 & \ldots & 0 & \ldots & 0 & 0 & \ldots & 0 & -1 & \ldots & \iX-1 \\
\end{array}%
\right),
$$
which is formed from $k$ blocks of sizes $(d_i-1)\times d_i$ and one last block of size $\iX\times\iX$.

\begin{definition}
\label{definition: number of components}
Let $l$ be a number of integral points in the convex hull of rows of the matrix $M_{d_1,\ldots,d_k;\iX}$ as elements in $\ZZ^{N-k}$.
Then \emph{the number $k_{d_1,\ldots,d_k;\iX}$} is $l-1$.
\end{definition}

Integral points of $\nabla$ are the origin and the points on the boundary.
The fan of $\widetilde T_X^\vee$ is given by a triangulation, vertices of whose simplices are integral points of $\nabla$
lying on the boundary. This means that boundary divisors of $\widetilde T_X^\vee$ are in one-to one correspondence
with the integral boundary points. Moreover, exceptional divisors of the resolution of $B$ do not lie
in the fiber over infinity. This means that the number of components of the fiber over infinity for $Z$
is equal to $k_{d_1,\ldots,d_k;\iX}$, and combinatorially they form a triangulation of a sphere
(which means that vertices of simplices of the triangulation correspond to components of the fiber over infinity,
edges correspond to intersections of components, etc.). In particular, the fiber over infinity is reduced,
so the degeneration of the fibers to the fiber over infinity is a semistable reduction,
and a monodromy of the family around infinity is maximally unipotent.
The compactification $Z$ is a tame compactified Landau--Ginzburg model, see~\cite{KKP17}.

\begin{problem}
Find a formula for $k_{d_1,\ldots,d_k;\iX}$ in terms of $d_1,\ldots,d_k,\iX$.
\end{problem}

The author is grateful to V.\,Golyshev for helpful comments.

\end{document}